\magnification=\magstep1
\documentstyle{amsppt}
\TagsOnLeft \CenteredTagsOnSplits
\pagewidth{16 true cm}\pageheight{22 true cm}
\input epsf.sty

\topmatter
\title Anisotropic Young diagrams and Jack symmetric functions \endtitle
\author S.~Kerov \endauthor
\address
Steklov Math. Institute (POMI), Fontanka 27,
St.Petersburg, 191011, Russia
\endaddress
\email kerov\@pdmi.ras.ru \endemail
\thanks
Partially supported by the Federal Grant Program
``Integration'', No.~326.53, and by MSRI at Berkeley.
\endthanks
\abstract
We study the Young graph with edge multiplicities
$\varkappa_\alpha(\lambda,\Lambda)$ arising in a Pieri-type
formula $p_1(x)\,P_\lambda(x;\alpha)=
\sum_{\Lambda:\lambda\nearrow\Lambda}
\varkappa_\alpha(\lambda,\Lambda)\,P_\Lambda(x;\alpha)$ for Jack
symmetric polynomials $P_\lambda(x;\alpha)$ with a parameter
$\alpha$. Starting with $\dim_\alpha\varnothing=1$, we define
recurrently the numbers $\dim_\alpha\Lambda=
\sum\varkappa_\alpha(\lambda,\Lambda)\,\dim_\alpha\lambda$, and
we set $\varphi(\lambda)=
\prod_{b\in\lambda}\big(a(b)\alpha+l(b)+1\big)^{-1}$ (where
$a(b)$ and $l(b)$ are the arm- and leg-length of a box $b$).

New proofs are given for two known results. The first is the
$\alpha$-hook formula $\dim_\alpha\lambda=n!\,\alpha^n\,
\prod_{b\in\lambda}\big((a(b)+1)\alpha+l(b)\big)^{-1}$. Secondly,
we prove (for all $u,v\in\Bbb{C}$) the summation formula
$\sum_{\Lambda:\lambda\nearrow\Lambda}
(c_\alpha(b)+u)(c_\alpha(b)+v)
\varkappa_\alpha(\lambda,\Lambda)\,\varphi(\Lambda)=
(n\alpha+uv)\;\varphi(\lambda)$, where $c_\alpha(b)$ is the
$\alpha$-content of a new box $b=\Lambda\setminus\lambda$. If
$\alpha=1$, this identity implies the existence of an interesting
family of positive definite central functions on the infinite
symmetric group.

The approach is based on the interpretation of a Young diagram as
a pair of interlacing sequences, so that analytic techniques may
be used to solve combinatorial problems. We show that when dealing
with Jack polynomials $P_\lambda(x;\alpha)$, it makes sense to
consider {\it anisotropic Young diagrams} made of rectangular
boxes of size $1\times\alpha$.
\endabstract
\endtopmatter

\document

\subhead Introduction \endsubhead
The basic observation behind this paper is that the {\it analytic
notion of a pair of interlacing sequences is a natural
generalization of the combinatorial notion of Young diagram} (see
\cite{5} for more details and \cite{3}, \cite{4} for applications
of combinatorial methods to the analysis of interlacing
sequences). We show that simple analytic facts related to
interlacing sequences imply the hook formulae for dimensions of
(representations corresponding to) Young diagrams, and for
transition probabilities of the Plancherel measure of the
infinite symmetric group.  More generally, we obtain transition
probabilities for the family of $z$-measures introduced in
\cite{6} which forms a deformation of the Plancherel measure.
Both results admit a straightforward generalization in which the
branching of Schur functions is replaced with that of Jack
symmetric polynomials.

The analytic approach to the hook formula for ordinary dimensions
was the subject of the papers \cite{1}, \cite{8}, \cite{11}.
The point of this note is that the considerations of \cite{1}
provide, with no extra effort, similar facts for the dimensions
related to Jack symmetric functions. The only difference is that
we identify Young diagrams with interlacing sequences in another
way, taking into account the parameter $\alpha$ of Jack
polynomials. Roughly, the ordinary Young diagrams should be
dilated by the factor of $\alpha$ along the horizontal axis. One
can also say that the new ``diagrams'' are built with rectangle
``boxes'' of width $\alpha$ and unit height, instead of unit
square boxes in case of ordinary Young diagrams. We call them
{\it anisotropic Young diagrams}.

The paper is organized as follows. We start in Section 1 by
recalling the identification of a Young diagram with a pair of
interlacing sequences. In Section 2 we use partial fractions to
associate with a pair of interlacing sequences two discrete
probability distributions. As explained in Sections 3-4, these
distributions generalize the transition and co-transition
distributions
$$
p_\lambda(\Lambda) = \frac{\dim\Lambda}{|\Lambda|\dim\lambda},
\quad q_\Lambda(\lambda) = \frac{\dim\lambda}{\dim\Lambda};
\qquad \lambda\nearrow\Lambda
$$
defined in combinatorial terms. In Section 5 we compute the
moments of the transition distribution assocoated with a pair of
interlacing sequences and show that the first three moments
depend not on the diagram itself but only on its number of boxes.
In Section 6 we give a short analytic proof of the $\alpha$-hook
formula first found by Stanley in \cite{10}. In the final Section
7 we establish our second main result, Theorem 7.5. It
generalizes to the case of Jack symmetric polynomials the family
of central measures on the Young lattice (or, equivalently, the
family of spherical functions on the infinite symmetric group)
introduced in \cite{6}.

{\bf Acknowledgement.}
I appreciate very much indeed the numerous discussions I had with
G.~I.~Olshanskii on the relationship between combinatorics and
analysis of Young diagrams, and on the combinatorics of Jack
polynomials in particular. I am also indebted to S.~Fomin,
A.~Okounkov and A.~Vershik for their interest in this piece of
work. The paper was completed during my visit to MSRI, Berkeley.
It is my pleasure to thank the organizers of the Combinatorics
1997 Program for the invitation and financial support.

\subhead 1. Interlacing sequences and Young diagrams \endsubhead
The aim of this Section is to explain that Young diagrams can be
naturally regarded as integer interlacing sequences.

Given a finite set $\{x_1,\,\ldots,x_d\}$ of distinct real
numbers, we shall always enumerate its elements in the increasing
order, and identify the set with the corresponding sequence.

\smallskip\noindent
(1.1) {\bf Definition.}
Two sets (or increasing sequences) $y_1,\,\ldots,y_{d-1}$ and
$x_1,\,\ldots,x_{d-1},x_d$ are said to be {\it interlacing}, iff
$$
x_1 < y_1 < x_2 < \ldots < x_{d-1} < y_{d-1} < x_d.
\tag 1.2
$$
The number $c=\sum x_k-\sum y_k$ is called the {\it center} of
interlacing sequences.

\smallskip
With every pair of interlacing sequences (1.2) we associate a
piecewise linear continuous function $v=\omega(u)$, such that
$$
\aligned
 (i)\qquad \omega'(u) &= +1,\quad
\text{ if } x_k < u < y_k,\quad k = 1,\ldots,d-1;\\
(ii)\qquad \omega'(u) &= -1,\quad
\text{ if } y_k < u < x_{k+1},\quad k = 1,\ldots,d-1;\\
(iii)\qquad \omega\,(u) &= |u - c|,\quad
\text{ if }\, u < x_1\, \text{ or }\, u > x_d.
\endaligned
\tag 1.3
$$
One can easily see that such a function exists, and is uniquely
determined by the properties (1.3). In fact, it follows from
(iii), that $\omega(x_1)=c-x_1$ and $\omega(x_d)=x_d-c$. By (i)
and (ii) we know that $\omega(y_k)-\omega(x_k)=y_k-x_k$ and
$\omega(x_{k+1})-\omega(y_k)=y_k-x_{k+1}$. This implies that
$$
\aligned
\omega(x_k) &=
\sum_{i=1}^{k-1} (y_i-x_i) + \sum_{j=k}^{d-1} (x_{j+1}-y_j), \\
\omega(y_k) &=
\sum_{i=1}^k (y_i-x_i) + \sum_{j=k}^{d-1} (x_{j+1}-y_j).
\endaligned
\tag 1.4
$$

The region
$$
S_\omega = \{(u,v)\in\Bbb{R}^2:\; |u - c| \le v < \omega(u)\}
\tag 1.5
$$
between the graphs of functions $v=\omega(u)$ and $v=|u-c|$
resembles the shape of a Young diagram, see Fig.~1. We say that
$\omega$ is the {\it diagram} of interlacing sequences (1.2), and
that the region (1.5) is its {\it shape}. We define the {\it
area} of this shape as
$$
A(\omega) = \sum_{i<j} (y_i-x_i)\,(x_j-y_{j-1}).
\tag 1.6
$$

We denote by $\Bbb{I}_d$ the space of interlacing sequences (1.2)
with $2d-1$ entries, and by $\Bbb{I}=\bigcup\Bbb{I}_d$ the set of
all interlacing sequences. The topology of the space $\Bbb{I}$ is
that of uniform convergence of associated diagrams.

$$
\epsffile{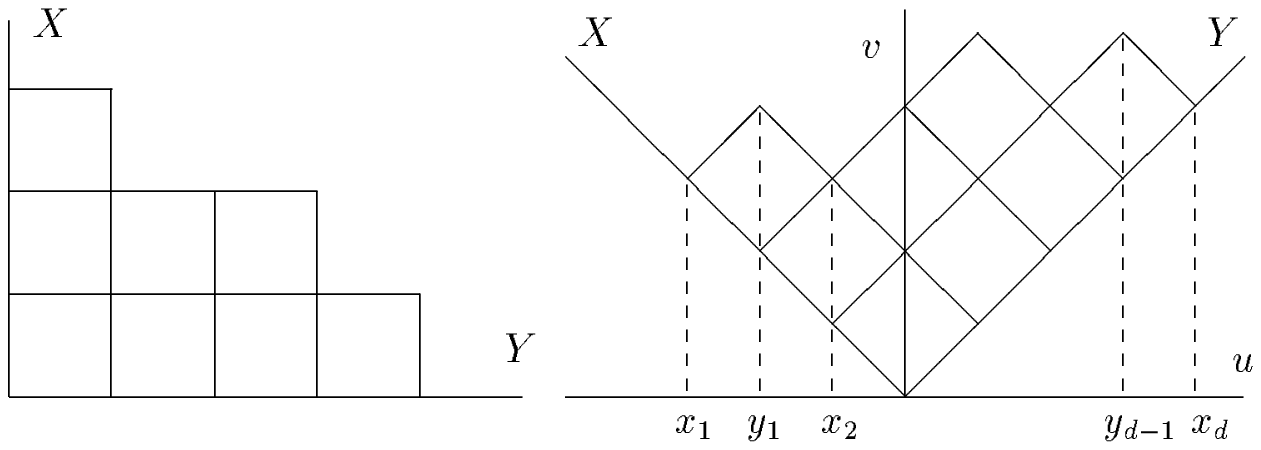}
$$
\centerline{Fig.~1.\quad
A Young diagram as a pair of interlacing sequences.}
\medskip

The shape of a true Young diagram
$\lambda=(\lambda_1,\,\ldots,\lambda_m)$ is uniquely determined
by the contents $y_1,\,\ldots,y_{d-1}$ of its corner boxes, and
the contents $x_1,\,\ldots,x_{d-1},x_d$ of corner boxes of the
compliment of $\lambda$ in $\Bbb{R}^2$. (Recall that the content
of a box $b=(i,j)$ on the crossing of the $i$-th row and $j$-th
column is defined as $c(b)=j-i$. It will also be convenient to say
that each point in the plane with the coordinates $(u,v)$ has the
{\it content} $c(u,v)=v-u$. Then the content of a box coincides
with that of its center.) The sequences of integers
$y_1,\,\ldots,y_{d-1}$ and $x_1,\,\ldots,x_{d-1},x_d$ interlace.
Also, the center of such a pair is always zero,
$c=\sum{x}_k-\sum{y}_k=0$.

In the opposite direction, every pair of interlacing sequences
(1.2) with integer entries and zero center represents a true Young
diagram. The number of rows $m$ in this diagram equals $-x_1$,
and $\lambda_1=x_d$ is the length of the first row. The number of
boxes $|\lambda|=A(\omega)$ equals the area of the associated
diagram $\omega$.

Therefore, the set $\Bbb{Y}$ of Young diagrams can be regarded as
the lattice of integer points in the space $\Bbb{I}_0$ of
interlacing sequences with zero center.

Conventionally, one draws a Young diagram $\lambda$ with
horizontal rows in decreasing order, so that the graph of the
associated piecewise linear function should be rotated by
$45^\circ$. The contents of the points of the horizontal parts on
the border of $\lambda$ form the intervals $(x_i,y_i)$, and those
of the vertical parts -- the intervals $(y_i,x_{i+1})$,
$i=1,\,\ldots,d-1$.

\subhead 2. Transition and co-transition distributions of an
interlacing sequence \endsubhead
In this Section we associate with every pair of interlacing
sequences two discrete probability distributions. In the
particular case of sequences corresponding to a Young diagram,
the distributions coincide with the transition and co-transition
distributions for the Plancherel measure of the infinite symmetric
group. All the facts that we recall here are well known, though
we use the combinatorial terminology related to the Young
lattice.

Given a pair (1.2) of sets $X=(x_1,\,\ldots,x_{d-1},x_d)$ and
$Y=(y_1,\,\ldots,y_{d-1})$, we denote by $P(u)=\prod(u-x_k)$ and
$Q(u)=\prod(u-y_k)$ the monic polynomials with corresponding
roots. Let us expand the rational fraction $R(u)=Q(u)/P(u)$ as a
sum of partial fractions,
$$
\frac{(u-y_1)\ldots(u-y_{d-1})}{(u-x_1)\ldots(u-x_{d-1})\,(u-x_d)}
= \sum_{k=1}^d \frac{\mu_k}{u-x_k}.
\tag 2.1
$$
Multiplying both sides by $u$ and taking the limit $u\to\infty$
we observe that
$$
\sum_{k=1}^d \mu_k = 1.
$$

\smallskip\noindent
(2.2) {\bf Lemma.} {\it
The following two properties of the fraction (2.1) are
equivalent:}
$$
x_1 < y_1 < x_2 < \ldots < x_{d-1} < y_{d-1} < x_d
\quad \Longleftrightarrow \quad
\mu_1 > 0,\; \ldots,\; \mu_d > 0.
$$

\demo{Proof}
By the residue formula,
$$
\mu_k = \frac{Q(x_k)}{P'(x_k)} =
\prod_{i<k} {x_k-y_i \over x_k-x_i}
\prod_{j>k} {x_k-y_{j-1} \over x_k-x_j}.
$$
Clearly, if the roots and poles of the fraction $R(u)$ interlace,
then all the coefficients $\mu_k$ are positive.

Assume now that $\mu_k>0$ for all $k=1,2,\,\ldots,d$. Then the
number of $y$-points to the right of the largest $x$-point $x_d$
is even. Moreover, it is zero. Otherwise, we could find two
neighboring $x$-points, not separated by a $y$-point, and the
corresponding coefficients $\mu$ would have opposite signs. A
similar argument proves that there is the same number of
$x$-points and $y$-points to the right of every $x_k$, so that
the sequences interlace.
\qed\enddemo

\smallskip\noindent
(2.3) {\bf Definition.}
To every pair of interlacing sequences
$x_1<y_1<\ldots<y_{d-1}<x_d$ there corresponds a discrete
probability distribution $\mu$ with the weights
$$
\mu_k = \prod_{i=1}^{k-1} {x_k-y_i \over x_k-x_i}
\prod_{j=k+1}^d {x_k-y_{j-1} \over x_k-x_j}
\tag 2.4
$$
at the points $x_k$, $k=1,2,\,\ldots,d$. We refer to $\mu$ as to
the {\it transition distribution} of the pair, and we call
the weights $\mu_k$ its {\it transition probabilities}.

Another distribution associated with a pair of interlacing
sequences arises from the decomposition
$$
\frac{(u-x_1)\ldots(u-x_{d-1})\,(u-x_d)}{(u-y_1)\ldots(u-y_{d-1})}
= u - c - \sum_{k=1}^{d-1} \frac{\nu_k}{u-y_k}.
\tag 2.5
$$

\smallskip\noindent
(2.6) {\bf Lemma.} {\it
The coefficient $c=c(\omega)$ in (2.5) coincides with the center
of the diagram $\omega$, and the sum of the coefficients
$\sum\nu_k=A(\omega)$ equals the area of its shape. The following
two conditions are equivalent:}
$$
x_1 < y_1 < x_2 < \ldots < x_{d-1} < y_{d-1} < x_d
\quad \Longleftrightarrow \quad
\nu_1 > 0,\; \ldots,\; \nu_{d-1} > 0.
\tag 2.7
$$

\demo{Proof}
In order to prove the first two claims of the Lemma we multiply
both sides of (2.5) by the polynomial $P(u)$ and compare a number
of first coefficients. Let $e_k(Y)$ and
$$
e_k(X) = \sum_{i_1<\ldots<i_k} x_{i_1} \ldots x_{i_k}
$$
denote the $k$-th elementary symmetric functions in the variables
$Y=(y_1,\,\ldots,y_{d-1})$ and $X=(x_1,x_2,\,\ldots,x_d)$
correspondingly. Up to the terms of degree $d-3$ and smaller,
$$
\gathered
u^d - e_1(X)u^{d-1} + e_2(X)u^{d-2} - \ldots =
-u^{d-2}\sum\nu_k + \ldots \\ +(u - c)\,
\big(u^{d-1} - e_1(Y)u^{d-2} + e_2(Y)u^{d-3} - \ldots\big),
\endgathered
$$
so that $c=e_1(X)-e_1(Y)$ and
$\sum\nu_k=-e_2(X)+e_1(X)e_1(Y)-e_1^2(Y)+e_2(Y)=A(\omega)$.

By the residue formula,
$$
\nu_k = - \frac{P(y_k)}{Q'(y_k)} =
- \prod_i (y_k-x_i) \prod_{j\ne k} (y_k-y_j)^{-1},
$$
If the sequences $X$ and $Y$ interlace, the number of $x$-points
to the right of each $y_k$ is one bigger than the number of
$y$-points to the right of $y_k$. It follows that the
coefficients $\nu_k$ are all positive.

Denote by $x(k)$, $y(k)$ the numbers of $x$-points and $y$-points
bigger than $y_k$. If the coefficient $\nu_k$ is positive, the
difference $x(k)-y(k)$ is odd. In fact, it can not be negative.
Otherwise, the biggest element would be $y_{d-1}$ or else we
could find two $y$-points not separated by an $x$-point. In both
cases some coefficient $\nu_m$, $m>k$, would be negative. The
assumption $x(k)-y(k)>1$ also implies a contradiction, for there
would be too many $y$-points to the left of $y_k$, and the sign
of $\nu_j$ for one of these points would be negative. Hence,
$x(k)=y(k)+1$ for all $k=1,\,\ldots,d-1$ which means that the
sequences $X$, $Y$ interlace.
\qed\enddemo

\smallskip\noindent
(2.8) {\bf Definition.}
Let $A=A(\omega)$ be the area (1.6) of the diagram $\omega$
associated with a pair of interlacing sequences (1.2). Denote by
$\nu$ the system of weights
$$
{\nu_k \over A} = \frac{(x_d-y_k)(y_k-x_1)}{A}\,
\prod_{i=1}^{k-1} \frac{y_k-x_{i+1}}{y_k-y_i}
\prod_{j=k+1}^{d-1} \frac{y_k-x_j}{y_k-y_j}
\tag 2.9
$$
at the points $y_k$, $k=1,\,\ldots,d-1$. By Lemma 2.6, $\nu$ is a
probability distribution. We call $\nu$ the {\it co-transition
distribution}, and the weights $\nu_k/A$ the {\it co-transition
probabilities} of the pair.

\subhead 3. Co-transition distribution of a Young diagram
\endsubhead
We show in this Section that the transition and co-transition
distributions introduced in Section 2 generalize the
corresponding combinatorial notions. We start by recalling the
definition of the Plancherel growth process (see \cite{2} and
\cite{12} for more details) and the related combinatorial
definitions of transition and co-transition probabilities.

We denote by $\Bbb{Y}_n$ the set of Young diagrams with $n$
boxes, and by $\Bbb{Y}=\bigcup\Bbb{Y}_n$ the lattice of all Young
diagrams ordered by inclusion. We write $\lambda\nearrow\Lambda$,
if the diagram $\Lambda$ covers $\lambda$ in $\Bbb{Y}$, i.e., if
$\Lambda=\lambda\bigcup{b_0}$ is a union of $\lambda$ and an
extra box $b_0$.

Let $\dim\lambda$ denote the number of standard Young tableaux of
shape $\lambda\in\Bbb{Y}$. The dimension function $\dim\lambda$
may be characterized by the initial condition
$\dim\varnothing=1$, and by the recurrence relation
$$
\dim\Lambda = \sum_{\lambda:\;\lambda\nearrow\Lambda} \dim\lambda
$$
(similar to that of the Pascal triangle) which can also be written as
$$
\sum_{\lambda:\;\lambda\nearrow\Lambda}\frac{\dim\lambda}{\dim\Lambda}=1.
\tag 3.1
$$
The ratios $q_\Lambda(\lambda)=\dim\lambda/\dim\Lambda$ are called {\it
co-transition probabilities} of the diagram $\Lambda$.

The celebrated {\it hook formula} for the dimension reads
$$
\dim\lambda = n!\; \prod_{b\in\lambda} h^{-1}(b),
\tag 3.2
$$
where $h(b)=(\lambda_i-j)+(\lambda'_j-i)+1$ is the hook length of
the box $b=(i,j)$ on the crossing of $i$-th row and $j$-th
column.

In order to prove (3.2), one can show that the right hand side
enjoys a recurrence relation similar to (3.1). We shall derive
this fact from Lemma 2.6. To this end, we first check that the
combinatorial definition of co-transition probabilities is
equivalent to the analytic one of Section 2.

\smallskip\noindent
(3.3) {\bf Lemma.} {\it
Denote by $x_1<y_1<\ldots<y_{d-1}<x_d$ the interlacing sequences
corresponding to a Young diagram $\Lambda$ with $A$ boxes, and
let $\lambda$ be the Young diagram obtained from $\Lambda$ by
removing a box with the content $y_k$. Let $f_\lambda$ denote the
right hand side of (3.2). Then}
$$
\frac{f_\lambda}{f_\Lambda} = \frac{(x_d-y_k)(y_k-x_1)}{A}\,
\prod_{i=1}^{k-1} \frac{y_k-x_{i+1}}{y_k-y_i}
\prod_{j=k+1}^{d-1} \frac{y_k-x_j}{y_k-y_j}.
\tag 3.4
$$

\demo{Proof}
By definition, the ratio in the left hand side can be written as
$$
{f_\lambda \over f_\Lambda} =
{1 \over |\lambda|} \prod_b \frac{h(b)+1}{h(b)},
$$
where the product runs over the boxes in the column of $\lambda$
above the new box $b_0=\Lambda\setminus\lambda$, and along the
row of $\lambda$ to the left of $b_0$. The hook lengths $h(b)$
are taken with respect to the smaller diagram $\lambda$. Consider
the intersection $I$ of the column with a block of equal rows of
Young diagram $\lambda$, see Fig.~2a. The block corresponds to a
vertical interval at the border of the diagram, and we denote by
$y_j$, $x_{j+1}$ the contents of its endpoints. If $j>k$, the
product along $I$ simplifies to
$$
\prod_{b\in I} \frac{h(b)+1}{h(b)} =
\frac{(y_j-y_k+1)}{(y_j-y_k)}
\frac{(y_j-y_k+2)}{(y_j-y_k+1)} \ldots
\frac{(x_{j+1}-y_k)}{(x_{j+1}-y_k-1)} =
\frac{(x_{j+1}-y_k)}{(y_j-y_k)}.
$$
If $j=k$, the smallest hook of $I$ equals $1$, and the product
reduces to $\prod_I(h+1)/h=(x_{k+1}-y_k)$.

In a similar way we consider the case where $I$ is the
intersection of the row containing the new box $b_0$ with a block
of equal columns of $\lambda$, see Fig.~2b. Such a block
corresponds to a horizontal interval on the border of $\lambda$,
and we denote by $x_i$, $y_i$ the contents of its endpoints. The
product along $I$ reduces in this case to
$$
\prod_{b\in I} \frac{h(b)+1}{h(b)} =
\frac{(y_k-y_i+1)}{(y_k-y_i)}
\frac{(y_k-y_i+2)}{(y_k-y_i+1)} \ldots
\frac{(y_k-x_i)}{(y_k-x_i-1)} =
\frac{(y_k-x_i)}{(y_k-y_i)}
$$
if $i<k$, and to $\prod_I(h+1)/h=(y_k-x_k)$ if $i=k$.
\qed\enddemo

The hook formula (3.2) readily follows from Lemmas (2.6) and
(3.3). In fact, these results imply that
$\sum{f}_\lambda/{f}_\Lambda=1$, hence the numbers $f_\lambda$
satisfy the same recurrence relations (3.1) as $dim\lambda$ and
the identity $\dim\lambda=f_\lambda$ follows.

$$
\epsffile{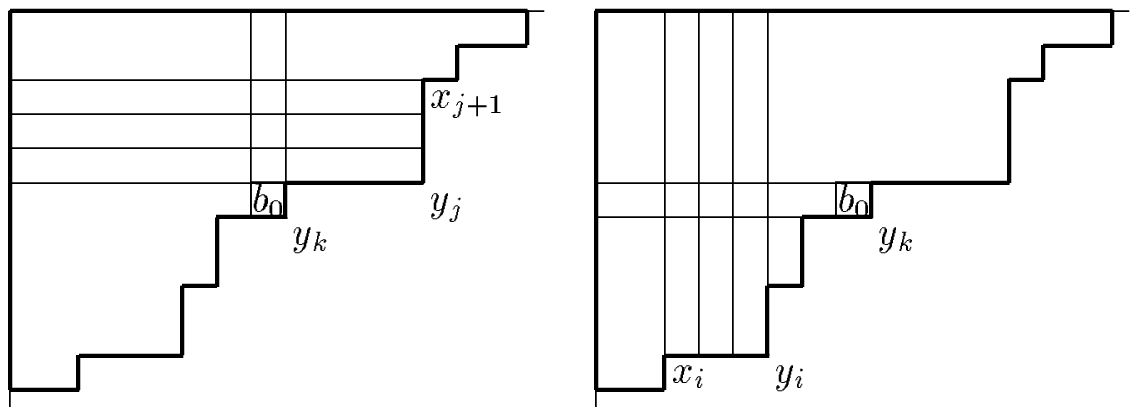}
$$
\centerline{Fig.~2a \hskip 50mm Fig.~2b}

\subhead 4. Plancherel transition distributions \endsubhead
Now we discuss the well known formula
$$
\sum_{\Lambda:\;\lambda\nearrow\Lambda} \dim\Lambda =
(|\lambda|+1)\dim\lambda
\tag 4.1
$$
which has a number of different proofs. For instance, one can
compute in two different ways the dimension of the representation
Ind\,$\pi_\lambda$ of the symmetric group $\frak{S}_{n+1}$
induced by the irreducible representation $\pi_\lambda$ (labeled
by a Young diagram $\lambda$) of the subgroup $\frak{S}_n$. A
bijective proof can be obtained from the basic properties of the
Robinson--Schensted insertion algorithm. Here we derive (4.1)
from the fact that the weights (2.4) associated with a pair of
interlacing sequences always form a probability distribution.

The formula (4.1) can be written as
$$
\sum_{\Lambda:\;\lambda\nearrow\Lambda}
{\dim\Lambda \over (n+1) \dim\lambda} = 1,
\tag 4.2
$$
so that the numbers
$\mu_\lambda(\Lambda)=\dim\Lambda/(n+1)\dim\lambda$, where
$\lambda\nearrow\Lambda$, can be taken as the weights of a
discrete probability distribution $\mu_\lambda$ associated with a
Young diagram $\lambda\in\Bbb{Y}_n$.

Consider the Markov chain on the Young lattice called the {\it
Plancherel growth process}. By definition, it starts at the empty
diagram $\varnothing$, and has the transition probabilities
$p(\lambda,\Lambda)=\dim\Lambda/(n+1)\dim\lambda$. One can easily
check that the probability to cross the $n$-th level $\Bbb{Y}_n$
of the Young lattice at a particular Young diagram $\lambda$
equals $M_n(\lambda)=\dim^2\lambda/n!$, the weight of this
diagram with respect to the Plancherel measure of the group
$\frak{S}_n$. In many ways, the chain can be regarded as the
Plancherel measure of the infinite symmetric group
$\frak{S}_\infty$. Following \cite{7}, we call $\mu_\lambda$
{\it Plancherel transition distribution}.

The Plancherel growth process is a {\it central} Markov chain,
meaning that the probability to reach a Young diagram $\lambda$
along a fixed path (i.e., Young tableaux) depends only on
$\lambda$, not on the choice of the path. In our case this
probability equals $\dim\lambda/n!$.

Recall that the {\it co-transition probability}
$q(\lambda,\Lambda)$ of a Markov chain on the Young lattice is
defined as the conditional probability to pass through a vertex
$\lambda$, assuming that the next vertex is $\Lambda$. For every
central Markov chain the co-transition probability is
$q(\lambda,\Lambda)=\dim\lambda/\dim\Lambda$. In fact, this
latter formula characterizes central Markov chains.

\smallskip\noindent
(4.3) {\bf Lemma.} {\it
Denote by $x_1<y_1<x_2<\ldots<y_{d-1}<x_d$ the interlacing
sequences associated with a Young diagram $\lambda\in\Bbb{Y}_n$.
Assume that a diagram $\Lambda$, where $\lambda\nearrow\Lambda$,
is obtained from $\lambda$ by attaching a box $b_0$ with the
content $x_k$. Then
$$
{\dim\Lambda \over (n+1) \dim\lambda} =
\prod_{i=1}^{k-1} {x_k-y_i \over x_k-x_i}
\prod_{j=k+1}^d {x_k-y_{j-1} \over x_k-x_j},
\tag 4.4
$$
i.e., combinatorial and analytic definitions of transition
probabilities are equivalent.}

\demo{Proof}
By the hook formula (3.2), the left hand side can be written as
$$
{\dim\Lambda \over (n+1) \dim\lambda} =
\prod_b \frac{h(b)}{h(b)+1},
$$
where the product runs over the boxes of $\lambda$ in the row, as
well as in the column, containing the new box $b_0$. All the hook
lengths $h(b)$ are taken with respect to the diagram $\lambda$.

The remaining argument is quite similar to that of the proof of
Lemma 3.3. Consider, for instance, the intersection $I$ of the
column through the box $b_0$ with a block of equal rows of
$\lambda$. Let $y_{j-1}$, $x_j$ be the contents of the endpoints
of the vertical interval on the crossing of these rows with the
graph of the diagram $\lambda$. Then
$$
\prod_{b\in I} \frac{h(b)}{h(b)+1} =
\frac{(y_{j-1}-x_k)}{(y_{j-1}-x_k+1)}
\frac{(y_{j-1}-x_k+1)}{(y_{j-1}-x_k+2)} \ldots
\frac{(x_j-x_k-1)}{(x_j-x_k)} =
\frac{(x_k-y_{j-1})}{(x_k-x_j)}.
$$
In a similar way we deal with the product along the part $I$ of
the row through $b_0$, corresponding to a horizontal interval
with the endpoints $x_i$, $y_i$ on the border of $\lambda$. It is
equal to
$$
\prod_{b\in I} \frac{h(b)}{h(b)+1} =
\frac{(x_k-y_i)}{(x_k-y_i+1)}
\frac{(x_k-y_i+1)}{(x_k-y_i+2)} \ldots
\frac{(x_k-x_i-1)}{(x_k-x_i)} =
\frac{(x_k-y_i)}{(x_k-x_i)},
$$
and the Lemma follows.
\qed\enddemo

$$
\epsffile{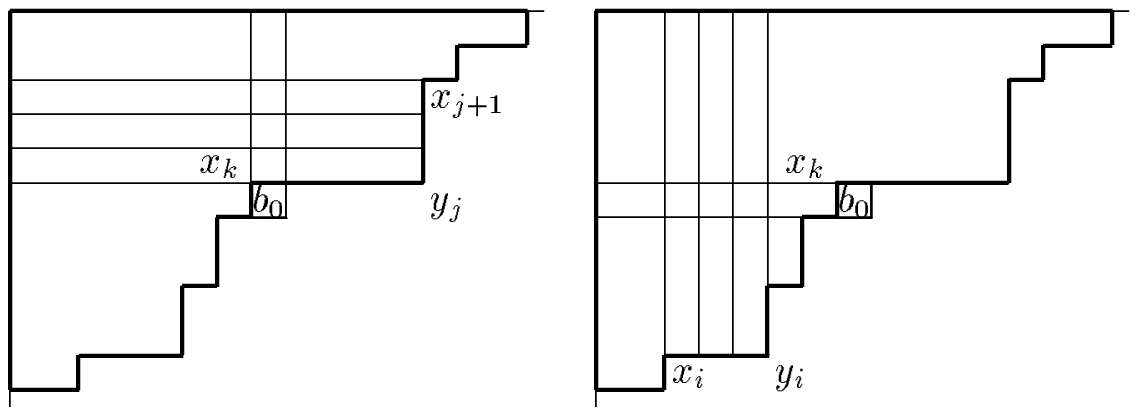}
$$
\centerline{Fig.~3a \hskip 50mm Fig.~3b}

\subhead 5. Transition distributions of $z$-measures \endsubhead
In this Section we find the moments of the transition
distribution (2.4) of a pair of interlacing measures. Then we use
these moments to define a family of central Markov chains on the
Young lattice, introduced in \cite{6}.

Let $\mu$ denote the transition distribution (2.4) of a pair
$x_1<y_1<x_2<\ldots<x_d$ of interlacing sequences, and
let
$$
h_m = \sum_{k=1}^d x_k^m\, \mu_k.
\tag 5.1
$$
be the $m$-th moment of $\mu$. We consider the generating function
$$
R(u) = \sum_{m=0}^\infty {h_m \over u^{m+1}}.
\tag 5.2
$$

\smallskip\noindent
(5.3) {\bf Lemma.} {\it
The moment generating function (5.2) equals}
$$
\sum_{m=0}^\infty {h_m \over u^{m+1}} =
\frac{(u-y_1)\ldots(u-y_{d-1})}
{(u-x_1)\ldots(u-x_{d-1})\,(u-x_d)}.
\tag 5.4
$$

\demo{Proof}
Since $1/(u-x_k)=\sum_{m\ge0}x_k^m/u^{m+1}$, we obtain that
$$
R(u) = \sum_{k=1}^d \frac{\mu_k}{u-x_k},
$$
and the claim follows from the identity (2.1).
\qed\enddemo

\smallskip\noindent
(5.5) {\bf Corollary.} {\it
Let $c$ and $A$ denote the center and the area of the diagram
associated with interlacing sequences $(X,Y)$. Then the mean
value of their transition distribution is $h_1=c$, and its
variance is $h_2-h_1^2=A$.}

\demo{Proof}
Using the notations of Section 2 for the elementary symmetric
functions, we derive from (5.4) that $h_1=e_1(X)-e_1(Y)=c$ and
$h_2=e_1^2(X)-e_2(X)-e_1(X)e_1(Y)+e_2(Y)$. Hence,
$h_2-h_1^2=-e_2(X)+e_1(X)e_1(Y)+e_2(Y)-e_1^2(Y)=A$.
\qed\enddemo

\smallskip\noindent
(5.6) {\bf Corollary.} {\it
Given two complex parameters $u$, $v$, consider the numbers
$$
p_{u,v}(\lambda,\Lambda) = \frac{(c(b)+u)(c(b)+v)}{n+uv}
\frac{\dim\Lambda}{(n+1)\dim\lambda}\;,
\tag 5.7
$$
where $\lambda\in\Bbb{Y}_n$, $\lambda\nearrow\Lambda$ and
$b=\Lambda\setminus\lambda$. Then}
$$
\sum_{\Lambda:\;\lambda\nearrow\Lambda}
p_{u,v}(\lambda,\Lambda) = 1.
$$

If $v=\bar{u}$, the numbers (5.7) are positive, and may be taken
as transition probabilities of a Markov chain on the Young
lattice. All these chains are central. In the limit $u\to\infty$
we obtain the Plancherel growth process, hence the family may be
considered as a deformation of the Plancherel measure of the
infinite symmetric group $\frak{S}_\infty$.

\subhead 6. Co-transition probabilities determined by the Jack
symmetric functions \endsubhead \break
Following \cite{9}, Chapter VI.10, we denote by
$P_\lambda(x)=P_\lambda(x_1,x_2,\ldots;\alpha)$;
$\lambda\in\Bbb{Y}$, the family of Jack symmetric functions with
the parameter $\alpha$. For a fixed value of $\alpha$ the
functions $P_\lambda$ form a linear basis in the symmetric
function algebra.

Denote by $p_1(x)=x_1+x_2+\ldots$ the sum of variables, and
consider the decomposition of the product $p_1\,P_\lambda$ in the
basis $P_\Lambda$. It is known (\cite{9}, VI.6.24(iv) and V.10.10)
that
$$
p_1(x)\, P_\lambda(x;\alpha) = \sum_{\Lambda:\lambda\nearrow\Lambda}
\varkappa_\alpha(\lambda,\Lambda)\, P_\Lambda(x;\alpha),
\tag 6.1
$$
where the multiplicities $\varkappa_\alpha(\lambda,\Lambda)$ are
given by an explicit formula
$$
\varkappa_\alpha(\lambda,\Lambda) = \prod_{b\in ver} \frac
{\big(a(b)\alpha+l(b)+2\big)\,\big(a(b)\alpha+l(b)+1\big)}
{\big((a(b)+1)\alpha+l(b)+1\big)\,\big(a(b)\alpha+l(b)+1\big)}\;.
\tag 6.2
$$
Here $b$ runs over all boxes in the $j$-th column $ver$ of the
diagram $\lambda$, provided that the new box
$b_0=\Lambda\setminus\lambda$ belongs to the $j$-th column of
$\Lambda$. The number $a(i,j)=\lambda_i-j$ is called the {\it arm
length}, and $l(i,j)=\lambda_j'-i$ is the {\it leg length} of a
box $b=(i,j)$ in the diagram $\lambda$.

We define a generalization $\dim_\alpha\lambda$ of the dimension
function $\dim\lambda$ by the recurrent formula
$$
\dim_\alpha\Lambda = \sum_{\lambda:\lambda\nearrow\Lambda}
\varkappa_\alpha(\lambda,\Lambda)\,\dim_\alpha\lambda,
\tag 6.3
$$
along with the initial condition $\dim_\alpha(\varnothing)=1$. If
$\alpha=1$, the multiplicities are trivial,
$\varkappa_1(\lambda,\Lambda)\equiv1$, and
$\dim_1\lambda=\dim\lambda$ is the ordinary dimension function.

The following hook-type formula for the dimension
$\dim_\alpha\lambda$, $|\lambda|=n$, was found by R.~Stanley
\cite{10}:
$$
\dim_\alpha\lambda = n!\,\alpha^n \prod_{b\in\lambda}
\Big(\big(a(b)+1\big)\alpha+l(b)\Big)^{-1}.
\tag 6.4
$$
We shall derive (6.4) as a particular case of Lemma 2.6.

To this end, we denote by $D_\alpha$ the dilation transform
$D_\alpha(u,v)=(u,v\alpha)$, $(u,v)\in\Bbb{R}^2$, and we consider
the image $\Lambda^\alpha=D_\alpha(\Lambda)$ of a Young diagram
$\Lambda$ upon this dilation. Note that the $\alpha$-{\it
content} $c_\alpha(u,v)=v\alpha-u$ of a point $(u,v)$ coincides
with the ordinary content $c(u,v\alpha)$ of the corresponding
point $D_\alpha(u,v)$. We denote by
$$
x_1(\alpha)<y_1(\alpha)<x_2(\alpha)<\ldots<
x_{d-1}(\alpha)<y_{d-1}(\alpha)<x_d(\alpha)
\tag 6.5
$$
the pair of interlacing sequences corresponding to the dilated
shape $\Lambda^\alpha$. In other words, (6.5) is the sequence of
$\alpha$-contents of corner points of the initial Young diagram
$\Lambda$. The center of the pair (6.5) is trivial,
$c(\alpha)=0$, and its area is $A(\alpha)=\alpha\,|\Lambda|$.

\smallskip\noindent
(6.6) {\bf Definition.}
Let $\Lambda$ be a Young diagram, $\lambda\nearrow\Lambda$, and
assume that $\alpha>0$. By (6.3), the numbers
$\varkappa_\alpha(\lambda,\Lambda)\dim_\alpha\lambda/\dim_\alpha\Lambda$
can be regarded as probabilities. We call them $\alpha$-{\it
co-transition} probabilities of the diagram $\Lambda$.

\smallskip\noindent
(6.7) {\bf Theorem.} {\it
Assume that a Young diagram $\lambda$ is obtained from
$\Lambda\in\Bbb{Y}_n$ by erasing a box $b_0$ with the
$\alpha$-content $y_k(\alpha)$ (more precisely, this is the
$\alpha$-content of the south-east corner point of the box). Then
$$
\aligned
{\varkappa_\alpha(\lambda,\Lambda)\,\dim_\alpha\lambda \over
\dim_\alpha\Lambda} &=
\frac{(x_d(\alpha)-y_k(\alpha))(y_k(\alpha)-x_1(\alpha))}
{A(\alpha)} \times \\
&\times \prod_{i=1}^{k-1} \frac{y_k(\alpha)-x_{i+1}(\alpha)}
{y_k(\alpha)-y_i(\alpha)}
\prod_{j=k+1}^{d-1} \frac{y_k(\alpha)-x_j(\alpha)}
{y_k(\alpha)-y_j(\alpha)}\,,
\endaligned
\tag 6.8
$$
i.e., the $\alpha$-co-transition probabilities of a Young diagram
$\Lambda$ coincide with the corresponding co-transition
probabilities of the associated pair of interlacing sequences
(6.5).}

\demo{Proof}
The proof is very similar to that of Lemma 3.3. By the
$\alpha$-hook formula (6.4),
$$
\frac{\dim_\alpha\lambda}{\dim_\alpha\Lambda} = {1 \over n}
\prod_{b\in ver} \frac{\big(a(b)+1\big)\alpha+l(b)+1}
{\big(a(b)+1\big)\alpha+l(b)}
\prod_{b\in hor} \frac{\big(a(b)+2\big)\alpha+l(b)}
{\big(a(b)+1\big)\alpha+l(b)},
$$
where $ver$ denotes the column of $\lambda$ above the new box
$b_0=\Lambda\setminus\lambda$, $hor$ is the row of $\lambda$ to
the left of $b_0$, and all arm- and leg-lengths are taken with
respect to the smaller diagram $\lambda$. We also took into
account the $\alpha$-hook
$h_\alpha(b_0)=\big((a(b)+1)\alpha+l(b)\big)=\alpha$ of the box
$b_0$.

Combining this with the multiplicity formula (6.2), we derive
that
$$
\frac{\varkappa_\alpha(\lambda,\Lambda)\,\dim_\alpha\lambda}
{\dim_\alpha\Lambda} = {1 \over n}
\prod_{b\in ver} \frac{a(b)\alpha+l(b)+2}
{a(b)\alpha+l(b)+1}
\prod_{b\in hor} \frac{\big(a(b)+2\big)\alpha+l(b)}
{\big(a(b)+1\big)\alpha+l(b)}.
\tag 6.9
$$

As we already did it in the proof of Lemma 3.3, consider the
intersection $I$ of $ver$ with a block of equal rows of the
diagram $\lambda$. Observe that the arm lengths along $I$ are
constant, and that the leg lengths increase by $1$. Assuming that
the rows of $I$ correspond to the interval
$\big(y_j(\alpha),x_{j+1}(\alpha)\big)$, we obtain that
$$
\prod_I \frac{a(b)\alpha+l(b)+2}{a(b)\alpha+l(b)+1} =
\frac{x_{j+1}(\alpha)-y_k(\alpha)}{y_j(\alpha)-y_k(\alpha)},
\qquad \text{\rm if } j > k,
$$
and
$$
\prod_I \frac{a(b)\alpha+l(b)+2}{a(b)\alpha+l(b)+1} =
x_{k+1}(\alpha)-y_k(\alpha),
\qquad \text{\rm if } j = k.
$$

For a part $I$ of the row $hor$ corresponding to the interval
$(x_i(\alpha),y_i(\alpha))$ we see that
$$
\prod_I \frac{\big(a(b)+2\big)\alpha+l(b)}
{\big(a(b)+1\big)\alpha+l(b)} =
\frac{y_k(\alpha)-x_i(\alpha)}{y_k(\alpha)-y_i(\alpha)},
\qquad \text{\rm if } i < k,
$$
and
$$
\prod_I \frac{\big(a(b)+2\big)\alpha+l(b)}
{\big(a(b)+1\big)\alpha+l(b)} =
\frac{y_k(\alpha)-x_k(\alpha)}{\alpha},
\qquad \text{\rm if } i = k.
$$

Since $A(\alpha)=n\alpha$, (6.8) follows from (6.9) and the
Theorem is proved.
\qed\enddemo

In the course of the proof of the Theorem we employed the
$\alpha$-hook formula (6.4) as a {\it definition} of the
$\alpha$-dimension function. We did not use any connection with
the recurrence formula (6.3). On the contrary, the hook formula
(6.4) is actually a direct corollary of Theorem 6.7 and Lemma
2.6. In fact, Lemma 2.6 implies that the function defined by
(6.4) satisfies the recurrence relation (6.3), hence the two
definitions are equivalent.

\smallskip\noindent
(6.10) {\bf Corollary \rm (Stanley, \cite{10}).} {\it
There is a hook formula
$$
\dim_\alpha\lambda = n!\,\alpha^n \prod_{b\in\lambda}
\Big(\big(a(b)+1\big)\alpha+l(b)\Big)^{-1}
$$
for the function $\dim_\alpha\lambda$ defined by (6.3).}

\subhead 7. The Plancherel growth process for anisotropic Young
diagrams \endsubhead
In the previous Section we considered a system of multiplicities
$\varkappa_\alpha(\lambda,\Lambda)$ for the edges of the Young
graph, depending on a parameter $\alpha>0$. The multiplicities
determine via the formula (6.3) a system of co-transition
probabilities. By Theorem 6.7, these probabilities
$$
q_\alpha(\lambda,\Lambda) =
{\varkappa_\alpha(\lambda,\Lambda)\,\dim_\alpha\lambda \over
\dim_\alpha\Lambda}, \qquad \lambda:\;\lambda\nearrow\Lambda
\tag 7.1
$$
also arise from the general Definition 2.8, if one replaces every
Young diagram $\Lambda$ by its dilated version $\Lambda^\alpha$.
Presently, we want to describe a distinguished Markov chain on
the Young lattice, which is central with respect to the
co-transition probabilities (7.1). It will be a generalization of
the Plancherel growth process for dilated (or {\it anisotropic})
Young diagrams.

\smallskip\noindent
(7.2) {\bf Lemma.} {\it
Consider a function $\varphi:\;\Bbb{Y}\to\Bbb{R}$ defined as
$$
\varphi(\lambda) = \prod_{b\in\lambda}
\big(a(b)\alpha+l(b)+1\big)^{-1}.
\tag 7.3
$$
Then
$$
\frac{\varkappa_\alpha(\lambda,\Lambda)\,\varphi(\Lambda)}
{\varphi(\lambda)} =
\prod_{i=1}^{k-1}
\frac{x_k(\alpha)-y_i(\alpha)}{x_k(\alpha)-x_i(\alpha)}
  \prod_{j=k+1}^d
\frac{x_k(\alpha)-y_{j-1}(\alpha)}{x_k(\alpha)-x_j(\alpha)},
\tag 7.4
$$
for every Young diagram $\lambda$.}

\demo{Proof}
It follows from the multiplicity formula (6.2) that
$$
\frac{\varkappa_\alpha(\lambda,\Lambda)\,\varphi(\Lambda)}
{\varphi(\lambda)} = \prod_{b\in ver}
\frac{\big(a(b)+1\big)\alpha+l(b)}{\big(a(b)+1\big)\alpha+l(b)+1}
\prod_{b\in hor}
\frac{a(b)\alpha+l(b)+1}{\big(a(b)+1\big)\alpha+l(b)+1},
$$
where $hor$ and $ver$ were defined in the proof of Theorem 6.7.

Let us denote by
$$
x_1(\alpha)<y_1(\alpha)<x_2(\alpha)<\ldots<
x_{d-1}(\alpha)<y_{d-1}(\alpha)<x_d(\alpha)
$$
the interlacing sequences corresponding to the dilated shape
$\lambda^\alpha$, and assume that the $\alpha$-content of (the
north-west corner point of) the new box
$b_0=\Lambda\setminus\lambda$ is $x_k(\alpha)$.

Proceeding as in the proof of Lemma 4.3, we split the column
$ver$ into blocks of boxes with equal arm lengths. The
$\alpha$-hooks $\big(a(b)+1\big)\alpha+l(b)$ in such a block $I$
vary from $y_{j-1}(\alpha)-x_k(\alpha)$ to
$x_j(\alpha)-x_k(\alpha)-1$, for some $j>k$. Hence,
$$
\prod_{b\in I}
\frac{\big(a(b)+1\big)\alpha+l(b)}{\big(a(b)+1\big)\alpha+l(b)+1}
= \frac{x_k(\alpha)-y_{j-1}(\alpha)}{x_k(\alpha)-x_j(\alpha)},
$$
and
$$
\prod_{b\in ver}
\frac{\big(a(b)+1\big)\alpha+l(b)}{\big(a(b)+1\big)\alpha+l(b)+1}
= \prod_{j=k+1}^d
\frac{x_k(\alpha)-y_{j-1}(\alpha)}{x_k(\alpha)-x_j(\alpha)}.
$$
Likewise,
$$
\prod_{b\in hor}
\frac{a(b)\alpha+l(b)+1}{\big(a(b)+1\big)\alpha+l(b)+1}
= \prod_{i=1}^{k-1}
\frac{x_k(\alpha)-y_i(\alpha)}{x_k(\alpha)-x_i(\alpha)},
$$
and the formula (7.4) follows.
\qed\enddemo

\smallskip\noindent
(7.5) {\bf Theorem.} {\it
The numbers $p_\alpha(\lambda,\Lambda)=
\varkappa_\alpha(\lambda,\Lambda)\,\varphi(\Lambda)/
\varphi(\lambda)$, $\lambda\nearrow\Lambda$, form a system of
transition probabilities for a Markov chain on the Young graph,
central with respect to the edge multiplicities (6.2).}

\demo{Proof}
By (7.4) and Lemma 2.2, the numbers $p_\alpha(\lambda,\Lambda)$
form, for every fixed Young diagram $\lambda$, a probability
distribution. The probability to reach a given Young diagram
$\Lambda\in\Bbb{Y}_n$ by a particular Young tableaux
$t=(\varnothing=\lambda_0\subset\lambda_1\subset\,
\ldots\subset\lambda_n=\Lambda)$ is easily seen to be
$$
P_\alpha(t) = \varphi(\Lambda)\,
\prod_{k=1}^n \varkappa_\alpha(\lambda_{k-1},\lambda_k),
$$
hence the co-transition probabilities are $q(\lambda,\Lambda)=
\varkappa_\alpha(\lambda,\Lambda)\,\dim_\alpha\lambda/\dim_\alpha\Lambda$,
for every $\Lambda$, and the chain is central.
\qed\enddemo

In complete analogy with Corollary 5.6 we derive from (7.4) and
Corollary 5.5 that the numbers
$$
p_{u,v}(\lambda,\Lambda) =
\frac{(c_\alpha(b)+u)(c_\alpha(b)+v)}{n\alpha+uv}\;
\frac{\varkappa_\alpha(\lambda,\Lambda)\,\varphi(\Lambda)}
{\varphi(\lambda)}
\tag 7.6
$$
sum up to unity. If all of these numbers are positive (e.g., if
$v=\bar{u}$) we can think of a Markov chain on the Young lattice
with transition probabilities (7.6).

\smallskip\noindent
(7.7) {\bf Corollary.} {\it
The Markov chain with the transition probabilities (7.6) is
central with respect to co-transition probabilities (7.1).}

\demo{Proof}
For every Young tableaux
$\varnothing=\lambda_0\subset\lambda_1\subset\,
\ldots\subset\lambda_n=\Lambda$, the product
$$
\prod_{k=1}^n \frac
{p_{u,v}(\lambda_{k-1},\lambda_k)}
{\varkappa_\alpha(\lambda_{k-1},\lambda_k)} =
{\varphi(\Lambda) \over uv(uv+\alpha)\ldots(uv+(n-1)\alpha)}
\prod_{b\in\lambda} (c_\alpha(b)+u)(c_\alpha(b)+v)
$$
depends on the final Young diagram $\Lambda$ only. This means
that the Markov chain is central.
\qed\enddemo


\enddocument
\end